\documentclass[12pt,a4paper]{amsart}
\theoremstyle{plain}

\usepackage{enumerate, amssymb}

\advance\hoffset-5mm \advance\textwidth40mm
\input{diagrams.tex}
\diagramstyle[scriptlabels,height=8mm,width=8mm]

\def\bdi{\begin{diagram}}
\def\edi{\end{diagram}}

\theoremstyle{plain}

\newtheorem{thm}{Theorem}[section]
\newtheorem{cor}[thm]{Corollary}
\newtheorem{lem}[thm]{Lemma}
\newtheorem{prop}[thm]{Proposition}
\theoremstyle{definition}
\newtheorem{defi}[thm]{Definition}
\newtheorem{defis}[thm]{Definitions}
\newtheorem{conj}[thm]{Conjecture}
\newtheorem{conv}[thm]{Convention}
\newtheorem{nota}[thm]{Notation}
\newtheorem{rem}[thm]{Remark}
\newtheorem{rems}[thm]{Remarks}
\newtheorem{exa}[thm]{Example}
\newtheorem{exas}[thm]{Examples}
\newtheorem{prob}[thm]{Problem}
\newtheorem{probs}[thm]{Problems}
\newtheorem{ques}[thm]{Question}
\newtheorem{sit}[thm]{}

\renewcommand{\epsilon}{\varepsilon}

\def\and{\quad\mbox{and}\quad}

\newcommand{\C}{\ensuremath{\mathbb{C}}}

\newcommand{\Z}{\ensuremath{\mathbb{Z}}}

\renewcommand{\rho}{\varrho}

\def\bals#1\eals{\begin{align*}#1\end{align*}}
\def\bal#1\eal{\begin{align}#1\end{align}}

\renewcommand{\phi}{\varphi}

\newcommand{\bnum}{\begin{enumerate}}
\newcommand{\enum}{\end{enumerate}}

\newcommand{\brem}{\begin{rem}}
\newcommand{\brems}{\begin{rems}}
\newcommand{\erem}{\end{rem}}
\newcommand{\erems}{\end{rems}}
\newcommand{\bprob}{\begin{prob}}
\newcommand{\eprob}{\end{prob}}
\newcommand{\bprobs}{\begin{probs}}
\newcommand{\eprobs}{\end{probs}}
\newcommand{\bques}{\begin{ques}}
\newcommand{\eques}{\end{ques}}
\newcommand{\bexa}{\begin{exa}}
\newcommand{\bexas}{\begin{exas}}
\newcommand{\eexa}{\end{exa}}
\newcommand{\eexas}{\end{exas}}
\newcommand{\bdefi}{\begin{defi}}
\newcommand{\edefi}{\end{defi}}
\newcommand{\bdefis}{\begin{defis}}
\newcommand{\edefis}{\end{defis}}
\newcommand{\bcor}{\begin{cor}}
\newcommand{\ecor}{\end{cor}}
\newcommand{\blem}{\begin{lem}}
\newcommand{\elem}{\end{lem}}
\newcommand{\bconv}{\begin{conv}}
\newcommand{\econv}{\end{conv}}
\newcommand{\bconj}{\begin{conj}}
\newcommand{\econj}{\end{conj}}
\newcommand{\bprop}{\begin{prop}}
\newcommand{\eprop}{\end{prop}}
\newcommand{\bthm}{\begin{thm}}
\newcommand{\ethm}{\end{thm}}
\newcommand{\bnota}{\begin{nota}}
\newcommand{\enota}{\end{nota}}
\newcommand{\bsit}{\begin{sit}}
\newcommand{\esit}{\end{sit}}
\newcommand{\be}{\begin{equation}}
\newcommand{\ee}{\end{equation}}
\newcommand{\bproof}{\begin{proof}}
\newcommand{\eproof}{\end{proof}}
\def\ba{\begin{array}}
\def\ea{\end{array}}

\thanks{ }

\begin{document}

\title[ One more proof of the Abhyankar-Moh-Suzuki Theorm]{ One more proof of the  Abhyankar-Moh-Suzuki  Theorm}

\author{S.\ Kaliman}
\address{Department of Mathematics,
University of Miami, Coral Gables, FL 33124, USA}
\email{kaliman@math.miami.edu}

\begin{abstract}   We extract the  Abhyankar-Moh-Suzuki  theorem from the Lin-Zaidenberg theorem. 
\end{abstract}





%

\begin{center}{\Large One more proof of the  Abhyankar-Moh-Suzuki  Theorem\\[3ex]}\end{center}

\begin{center}{Shulim Kaliman\\[3ex]}\end{center}

\begin{center}  ABSTRACT. We extract the  Abhyankar-Moh-Suzuki  theorem from the Lin-Zaidenberg theorem. \\[5ex]
\end{center}

\section*{Introduction}
Let $\Gamma_0$ be the zero locus of a primitive polynomial $p \in \C [x,y]$. Suppose that $\Gamma_0$ is a smooth irreducible simply connected curve
(i.e. $\Gamma_0 \simeq \C$). Then the Abhyankar-Moh-Suzuki  theorem  \cite{AbMo}, \cite{Su} states that

{\em  in a suitable polynomial coordinate system on $\C^2$ the curve $\Gamma_0$ is a coordinate
axis, i.e. $p$ may be viewed as a coordinate function. }

The Lin-Zaidenberg theorem \cite{LZ} states that if $\Gamma_0$ is not isomorphic but only homeomorphic to $\C$ as a topological space then

{\em  in a suitable polynomial coordinate system on $\C^2$ the polynomial $p(x,y)$ is of the form $x^k +y^l$ where $k, l\geq 2$ are relatively prime.}

These theorems are extremely important tools of the modern affine algebraic geometry and it is not a surprise that different methods of proving them were suggested after the original papers.
In the case of  the Abhyankar-Moh-Suzuki  theorem there is at least a dozen of proofs (e.g., see
 \cite{Mi},   \cite{Ru},   \cite{Ri}, \cite{Ne},  \cite{Kan},   \cite{ACL},   \cite{CO},  \cite{GM},  \cite{Gu}, \cite{Zo}, \cite{Ko}, \cite{ML} \cite{Pa}).
Some of these papers  contain actually simultaneous proofs of  the  Abhyankar-Moh-Suzuki and Lin-Zaidenberg  theorems (e.g., see \cite{Ne}, \cite{Ko}, \cite{Pa}) but
none of them used the beautiful technique 
developed in \cite{LZ} (and refined later in \cite{Za}).
At first glance this technique (based of the theory of Teichm\"uller spaces and hyperbolic analysis) does not imply the  Abhyankar-Moh-Suzuki  theorem but we show below that 
it is applicable and there is a simple reduction of the 
Abhyankar-Moh-Suzuki  theorem to the Lin-Zaidenberg one. 

{\em Acknowledgements.} In the first version of this text written about 30 years ago the author expressed his gratitude to V. Ya. Lin and M. Zaidenberg for useful discussions.
Though none of us remember these discussions now the gratitude still stands. The author is also indebted to the referee for his essential contribution to the quality of this paper.

\section{Reduction}

\bdefi\label{d1} Recall that every pair  $(k,l)$ of relatively prime 
natural numbers defines a weighted degree function $d$ on the ring of polynomials in $x$ and $y$ such that for a nonzero polynomial $r(x,y)=\sum_{i,j}a_{ij}x^iy^j$
one has  $d(r)=\max_{(i,j) \in I} (li+kj)$ where  $I=\{ (i,j) \in (\Z_{\geq 0})^2| a_{ij} \ne 0\}$. We let $\hat I = \{ (i,j) \in I| li+kj= d(p)\}$ and we call 
$\hat r (x,y) =\sum_{(i,j) \in \hat I}a_{ij}x^iy^j$ the leading quasi-homogeneous part of $r$ with respect to the grading induced
by $d$.
\edefi

\blem\label{l1} Let the zero locus of a primitive polynomial $p(x,y)$ be isomorphic to $\C$. Then there is an automorphism $\alpha$ of $\C^2$ 
 such that either  $q:= p \circ \alpha$ is linear or the leading quasi-homogeneous part of $q$ with respect to some weighted degree function $d$
is of the form $\hat q(x,y)=(x^k +y^l)^n$  where $k, l\geq 2$ are relatively prime.
\elem

\bproof 
For $d$ as in Definition \ref{d1} the fact that $k$ and $l$ are relatively prime implies that  
the leading quasi-homogeneous part  of any nonzero polynomial (and in particular $ p $) is automatically of the form $$ax^{n_0}\prod_{i=1}^m (c_ix^k+y^l)^{n_i}$$
where $n_i \geq 0$, $a \ne 0$, $c_1, \ldots , c_m \in \C$ are distinct, and $d(p)=n_0+kl \sum_{i =1}^mn_i$.
Furthermore, unless $p$ depends on one variable only (in which case it is automatically linear) the pair $(k,l)$ can be chosen so that $m\geq 1$.

However, in this case after the natural embedding
of $\C^2$ into the weighted projective plane, associated with the grading induced by $d$, one can see that the irreducible curve $\Gamma_0= p^{-1}(0)$
must have $m+1$ (resp. $m$) points at infinity for $n_0>0$ (resp. $n_0=0$) and, therefore, at least
the same number of punctures. Since $\Gamma_0 \simeq \C$ has one puncture we conclude that $$\hat p (x,y)=a(cx^k+y^l)^{n}$$ where $a$ and $c \ne 0$.
If either $k$ or $l$ (say $l$) is equal to 1 then applying a de Joinqu\`ere's transformation $(x,y) \to (x, y+cx^k)$ we can reduce the standard degree of $p$.
Thus, choosing $\alpha$ so that $q=p\circ \alpha$ has the minimal possible standard degree we see that unless $q$ is linear one can suppose that
$\hat q(x,y)=(x^k+y^l)^n$ with relatively prime $k, l \geq 2$.

\eproof

\bprop\label{l2} In Lemma \ref{l1} one has $n=1$, i.e.unless $q$ is linear the leading quasi-homogeneous part of $q$ is $\hat q(x,y)=x^k+y^l$
with relatively prime $k, l \geq 2$.

\eprop

\bproof Consider the surface $S$ given by $$z^{nkl}q(xz^{-l}, yz^{-k})=1$$ in $\C^3$ together with the family of curves $z|_S : S \to \C$. 
Then for $a\ne 0$ the curve $S_a= S\cap \{ z=a \}$ is isomorphic to
the fiber $\Gamma_c = \{ q(x,y)=c \}$ where $c=a^{-nkl}$. By \cite{LZ} all nonzero fibers of $q$ are pairwise isomorphic, which yields a pairwise isomorphism of curves $S_a$ for $a \ne 0$.
Note that $S_0$ is given by equation $\hat q(x,y)=1$ in the plane $\{ z=0 \} \subset \C^3$ and thus $S_0$ 
consists of $n$ disjoint components isomorphic to the curve $x^k+y^l=1$. However, in the case of isomorphic irreducible general fibers any degenerate fiber may contain
at most one hyperbolic component\footnote{More precisely, the smooth part of such a fiber contains at most one hyperbolic component. Furthermore, this component
is isomorphic to a quotient of the general fiber with
respect to a finite group action. } \cite[Theorem 5.3]{Za}.  Thus $n=1$. 
\eproof

\bcor\label{c0} If $q$ in Proposition \ref{l2} is not linear then general fibers of $q$ have negative Euler characteristic.
\ecor

\bproof Indeed, for relatively prime $k,l\geq 2$ the Riemann-Hurwitz formula implies that the  Euler characteristic of the curve $x^k+y^l=1$ (and, therefore, $S_0$)
is at most $-1$. The Suzuki formula for Euler characteristic \cite{Su} implies that the same is true for general curves $S_a$.
\eproof

\blem\label{l3} Let $\Delta_r$ be the disc $\{ x\in \C | \, |x| <r \}$ and $\pi : \C^2 \to \C$ be the projection to the $x$-axis. Then under the assumption of
Lemma \ref{l1} for every $\delta >0$ there exists $r >0$
such that $\pi$ maps the set $$R=\{ (x,y) \in \C^2 | \, \frac{\partial q(x,y)} {\partial y} =0\, \,  \& \, \, |q(x,y)|< \delta \}$$ into $\Delta_r$.

\elem

\bproof  By Proposition \ref{l2} $\hat q (x,y)=x^k +y^l$. Hence in the polynomial polyhedron $\{ (x,y) \in \C^2 | \, |q(x,y)|< \delta \}$ 
for $x$ with large absolute value one has $|x|^l \approx |y|^k$. Furthermore,
looking at the quasi-leading part of the polynomial $ \frac{\partial q(x,y)} {\partial y}$ one can see now that for such values of $x$ 
the absolute value of  $\frac{\partial q(x,y)} {\partial y}$ must be also large, which yields the desired conclusion.

\eproof

\blem\label{l4}
For $\delta$ and $r$ from Lemma \ref{l3} and any $c$ with $0< |c| <\delta$ the set $\Gamma_c \cap \pi^{-1}(\C \setminus \bar \Delta_r )$ is a disjoint union of components biholomorphic
to a punctured disc.

\elem

\bproof Note that  the solutions of the system $ \frac{\partial q(x,y)} {\partial y} =q(x,y)-c=0$ are ramification points of the finite morphism $\pi|_{\Gamma_c} : \Gamma_c \to \C$ 
that is the restriction of $\pi$. 
Hence by Lemma \ref{l3} $\Gamma_c \cap \pi^{-1} (\C \setminus \bar \Delta_r )$ is the union of components biholomorphic to a punctured disc.

\eproof

\bcor\label{l5} The Riemann surfaces $\Gamma_c$ and $\Gamma_c \cap \pi^{-1} (\Delta_r )$ are diffeomorphic. 

\ecor

\subsection{Proof of the  Abhyankar-Moh-Suzuki  theorem} Assume to the contrary that $q$ in  Lemma \ref{l1} and Proposition \ref{l2} is not linear. Then 
the leading quasi-homogeneous  part of $q$ is $x^k +y^l$ with $k,l \geq 2$.
Note that for large $r$ the set $\Gamma_0 \cap \pi^{-1}(\Delta_r )$ is biholomorphic to a disc.
The smoothness of $\Gamma_0$ implies that $\Gamma_c \cap \pi^{-1}(\Delta_r )$ is also biholomorphic to a disc for sufficiently small $|c|$. By Corollary \ref{l5}
$\Gamma_c$ is contractible in contradiction with Corollary \ref{c0}. Hence $q$
must be linear, which concludes the proof.  

\hspace{14.7cm} $\square$

\brem\label{r1}
The  Abhyankar-Moh-Suzuki  theorem is valid over any algebraically closed field $K$ of characteristic zero while the proof before is done over the field
of complex numbers. However, the Lefschetz Principle enables us to reduce the general case to the complex one. 

Indeed, the description of the polynomial $p$ and
an isomorphism between $p^{-1}(0)$ and a line involves a finite number of elements of $K$. These elements generate an algebraically closed subfield $K_0$ of $K$
which can be also embedded into $\C$. We can consider $p$ as a polynomial over $K_0$ and, consequently, over $\C$.  If it is known
already that every fiber $\Gamma_c = p^{-1}(c)$ of $p$ over $\C$ is isomorphic to a line then the same is true for the fibers of
$p$ over $K_0$ (and, therefore, over $K$) since the genus and the number of punctures of any affine curve survive the field extention $\C : K_0$.  
Hence one can extract the general form of the  Abhyankar-Moh-Suzuki  theorem from the fact that a polynomial with general fibers
isomorphic to a line is a variable in a suitable polynomial coordinate system \cite{Gut}.
\erem

\end{document}